\title{Unstable maps}

\author{Gareth A. Jones\\
School of Mathematics\\
University of Southampton\\
Southampton SO17  1BJ, UK\\
{\tt G.A.Jones@maths.soton.ac.uk}\\
}

\documentclass[12pt]{article}
\usepackage[latin1]{inputenc}
\usepackage{latexsym}
\usepackage{amsmath}
\usepackage{amssymb}
\usepackage{amsfonts}
\usepackage{color}
\usepackage{tikz}
\usepackage{enumerate}
\newtheorem{thm}{Theorem}[section]
\newtheorem{lemma}[thm]{Lemma}

\newcommand{\M}{\mathcal{M}}

\newcommand{\K}{\mathcal{K}}

\date{}
\begin{document} 

\maketitle

\begin{abstract}
A map which is non-orientable or has non-empty boundary has a canonical double cover which is orientable and has empty boundary. The map is called stable if every automorphism of this cover is a lift of an automorphism of the map. This note describes several infinite families of unstable maps, and relates them to similar phenomena for graphs, hypermaps and Klein surfaces.
\end{abstract}

\medskip

%05C10  Top. graph theory
%14H37 automorphisms (of curves)
%14H57 dessins d'enfants
% 20B25 finite automorphism groups
% 30F10 compact Riemann surfaces
% 30F50 Klein surfaces
 
 \medskip
 

\section{Introduction}

\begin{figure}[h!]
\begin{center}
\begin{tikzpicture}[scale=0.8, inner sep=0.8mm]

\draw [thick, dotted] (10,0) arc (0:360:3);
\node (a1) at (7,2) [shape=circle, fill=black] {};
\node (b1) at (8.9,0.618) [shape=circle, fill=black] {};
\node (c1) at (8.17,-1.616) [shape=circle, fill=black] {};
\node (d1) at (5.83,-1.616) [shape=circle, fill=black] {};
\node (e1) at (5.1,0.618) [shape=circle, fill=black] {};
\node (f1) at (7,0) [shape=circle, fill=black] {};
\draw [thick](a1) to (b1) to (c1) to (d1) to (e1) to (a1);
\draw [thick](a1) to (f1) to (b1);
\draw [thick](c1) to (f1) to (d1);
\draw [thick](e1) to (f1);

\draw [thick] (7.927,2.853) to (a1) to (6.033,2.853);
\draw [thick] (4.533,1.763) to (e1) to (4,0);
\draw [thick] (4.573,-1.763) to (d1) to (6.073,-2.853);
\draw [thick] (7.927,-2.853) to (c1) to (9.427,-1.763);
\draw [thick] (10,0) to (b1) to (9.467,1.763);

\draw (7.3,3.3) to (7,3) to (7.3,2.7);
\draw (6.7,-3.3) to (7,-3) to (6.7,-2.7);

\end{tikzpicture}

\end{center}
\caption{Regular embedding of $K_6$ in the real projective plane}
\label{projplK6}
\end{figure}
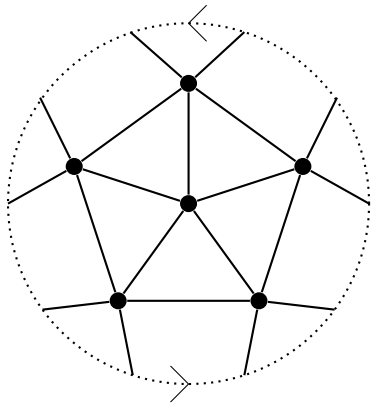

If a map $\M$ is non-orientable or has a non-empty boundary, then it has a {\em canonical orientable double cover\/} $\widetilde\M$, the unique double cover of $\mathcal M$ which is orientable and has an empty boundary. Then $\M$ is the quotient of $\widetilde\M$ by an orientation-reversing automorphism $a\in{\rm Aut}\,\widetilde\M$ of order $2$. In particular, if $\M$ has empty boundary then the covering $\widetilde\M\to\M$ is unbranched, with $\chi(\widetilde\M)=2\chi(\M)$. For example, the non-orientable regular embedding $\M$ of the complete graph $K_6$ in the real projective plane, shown in Figure~\ref{projplK6}, lifts to the icosahedral map $\widetilde\M$ on the sphere, and $\M$ is the quotient of $\widetilde\M$ by its antipodal automorphism $a$.

Each automorphism of $\M$ lifts to two automorphisms of $\widetilde\M$, differing by and commuting with $a$, so that ${\rm Aut}\,\widetilde\M$ has a subgroup ${\rm Aut}\,\M\times\langle a\rangle\cong{\rm Aut}\,\M\times C_2$. We will say that $\M$ is {\em stable\/} or {\em unstable\/} as this subgroup is the whole of ${\rm Aut}\,\widetilde\M$ or a proper subgroup. (This terminology is borrowed from Graph Theory, where it is used in a similar situation concerning bipartite double coverings of graphs \cite{LMS, MSZ, Wil08}.) While it appears that most maps are stable, the aim of this note is to give some examples of unstable maps, exhibiting various degrees of instability.

\medskip

The following example shows that if $\M$ has a non-empty boundary then $|{\rm Aut}\,\widetilde\M|$ can be arbitrarily large while $|{\rm Aut}\,\M|$ is bounded above.

\medskip

\noindent{\bf Example 1.} Let $\K$ be the regular spherical map $\{2,n\}$, in the notation of Coxeter and Moser~\cite{CM}: this has two vertices, joined to each other by $n$ edges, with ${\rm Aut}\,\K\cong D_n\times C_2$ of order $4n$. Let $a$ be the reflection of $\K$ in one of its edges, and let $\M=\K/\langle a\rangle$, a map on the closed disc. Then $\widetilde\M=\K$, and $|{\rm Aut}\,\M|=2$ or $4$ as $n$ is odd or even. For a more extreme example, take $\K$ to be a semi-star map on the sphere, with one vertex, one face and $n$ semi-edges, so that ${\rm Aut}\,\K\cong D_n$; if $\M=\K/\langle a\rangle$ where $a$ is a reflection of $\K$, then $\K=\widetilde\M$ and $|{\rm Aut}\,\M|=1$ or $2$ as $n$ is odd or even.

\medskip

One can easily find similar examples of higher genus, where $\M$ has a non-empty boundary.  However, examples with an empty boundary are not so obvious. In order to investigate these we will use a more algebraic approach.

%%%%%%%%%

\section{Algebraic map theory}

\begin{figure}[h!]
\begin{center}
\begin{tikzpicture}[scale=0.6, inner sep=0.8mm]

\node (c) at (0,0) [shape=circle, fill=black] {};
\node (d) at (8,0) [shape=circle, fill=black] {};
\draw [thick] (c) to (d);
\draw [thick] (c) to (1,-3);
\draw [thick] (c) to (1,3);
\draw [thick] (d) to (7,-3);
\draw [thick] (d) to (7,3);
\draw [thick] (c) to (-2.5,2.5);
\draw [thick] (c) to (-2.5,-2.5);
\draw [thick] (d) to (10.5,2.5);
\draw [thick] (d) to (10.5,-2.5);

\draw (c) to (1,0.5);
\draw (c) to (1,-0.5);
\draw (1,0.5) to (1,-0.5);
\draw (d) to (7,0.5);
\draw (d) to (7,-0.5);
\draw (7,0.5) to (7,-0.5);
\draw (c) to (0.8,0.8);
\draw (0.3,1) to (0.8,0.8);

\node at (-0.8,0) {$v$};
\node at (4,-0.4) {$e$};
\node at (4,2) {$f$};

\node at (1.5,0.6) {$\phi$};
\node at (6.2,0.6) {$\phi R_0$};
\node at (1.4,1.3) {$\phi R_1$};
\node at (1.9,-0.6) {$\phi R_2$};
\node at (5.9,-0.6) {$\phi R_0R_2$};

\end{tikzpicture}

\end{center}
\caption{Generators $R_i$ of $\Gamma$ acting on a flag $\phi=(v,e,f)$.} 
\label{flags}
\end{figure}
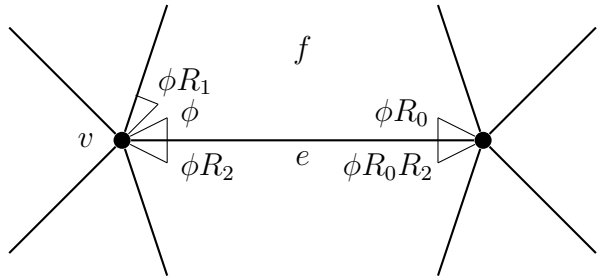

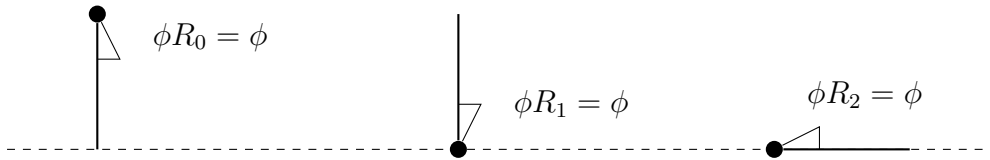
\begin{figure}[h!]
\begin{center}
\begin{tikzpicture}[scale=0.6, inner sep=0.8mm]

\draw [dashed] (-10,0) to (12,0);

\node (a) at (-8,3) [shape=circle, fill=black] {};
\draw [thick] (a) to (-8,0);
\draw (a) to (-7.5,2) to (-8,2);
\node at (-5.5,2.5) {$\phi R_0=\phi$};

\node (b) at (0,0) [shape=circle, fill=black] {};
\draw [thick] (b) to (0,3);
\draw (b) to (0.5,1) to (0,1);
\node at (2.5,1) {$\phi R_1=\phi$};

\node (c) at (7,0) [shape=circle, fill=black] {};
\draw [thick] (c) to (10,0);
\draw (c) to (8,0.5) to (8,0);
\node at (9,1.2) {$\phi R_2=\phi$};

\end{tikzpicture}

\end{center}
\caption{Flags fixed by $R_0, R_1$ and $R_2$.}
\label{fixedflags}
\end{figure}

As shown, for example, in~\cite{BS, Jon16, JT}, maps (always assumed to be connected) can be identified with transitive permutation representations of the group
\[\Gamma=\langle R_i\;(i=0,1,2)\mid R_i^2=(R_0R_2)^2=1\rangle\cong V_4*C_2.\]
Given any map $\M$, $\Gamma$ acts transitively on the set $\Phi$ of flags (incident vertex-edge-face triples) of $\M$, with $R_i$ changing the $i$-dimensional component of each flag (whenever possible) while preserving its other two components, as shown in Figure~\ref{flags}; flags incident with the boundary are fixed by some $R_i$, as shown in Figure~\ref{fixedflags} where the broken line represents the boundary of the map. Conversely, given a transitive permutation representation of $\Gamma$, one can define a map $\M$ by taking the vertices, edges and faces to be the orbits of the dihedral subgroups $\langle R_1, R_2\rangle$,  $\langle R_2, R_0\rangle$ and  $\langle R_0, R_1\rangle$, with incidence given by non-empty intersection.

In this situation the automorphisms of $\M$ are the permutations of $\Phi$ commuting with those induced by $\Gamma$; these form a group ${\rm Aut}\,\M\cong N_{\Gamma}(M)/M$, where $M$ is the stabiliser in $\Gamma$ of a flag, called a {\em map subgroup} for $\M$. The most symmetric maps are the {\em regular\/} maps, those for which ${\rm Aut}\,\M$ acts transitively on the flags; equivalently, $\Gamma$ induces a regular permutation group on $\Phi$, so that $M$ is a normal subgroup of $\Gamma$ and ${\rm Aut}\,\M\cong\Gamma/M$.

A map $\M$ is orientable and with empty boundary if and only if $M$ is contained in the even subgroup $\Gamma^+$ of index $2$ in $\Gamma$, consisting of the words of even length in the generators $R_i$. If $M\not\le\Gamma^+$, then the canonical double cover $\widetilde\M$ of $\M$ is the map with map subgroup $M^+:=M\cap\Gamma^+$, corresponding to the transitive action of $\Gamma$ on $\widetilde\Phi:=\Phi\times\{\pm 1\}$ given by $R_i:(\phi,\delta)\mapsto (\phi R_i,-\delta)$. Then $\M$ is the quotient of $\widetilde\M$ by the automorphism group of order $2$ corresponding to the inclusion $M/M^+\le N_{\Gamma}(M^+)/M^+\cong{\rm Aut}\,\widetilde\M$.

Since $M$ and $\Gamma^+$ are both normalised by $N_{\Gamma}(M)$, so is their intersection $M^+$, giving $N_{\Gamma}(M^+)\ge N_{\Gamma}(M)$. The group ${\rm Aut}\,\widetilde\M\cong N_{\Gamma}(M^+)/M^+$ therefore has a subgroup corresponding to
\[N_{\Gamma}(M)^+/M^+\times M/M^+\cong  N_{\Gamma}(M)/M\times M/M^+\cong {\rm Aut}\,\M\times C_2\]
obtained by lifting each automorphism of $\M$ to two automorphisms of $\widetilde\M$. Then $\M$ is stable if and only if this subgroup is the whole of ${\rm Aut}\,\widetilde\M$, or equivalently $N_{\Gamma}(M^+)=N_{\Gamma}(M)$.

%%%%%%%%%%

\section{Examples of unstable maps}

The following trivial observation indicates one area where we should not look for unstable maps.

\begin{lemma}
If $\mathcal M$ is regular then $\M$ is stable.
\end{lemma}

\noindent{\sl Proof.} If $\M$ is regular then $M$ is normal in $\Gamma$, so $N_{\Gamma}(M)=\Gamma$ and hence $N_{\Gamma}(M^+)=\Gamma$. Thus $N_{\Gamma}(M^+)=N_{\Gamma}(M)$, so $\M$ is stable. \hfill$\square$

\medskip

For instance, Example~1 illustrates this phenomenon. It is natural to ask whether weaker symmetry conditions also imply stability. A map $\M$ is {\em edge-transitive\/} if ${\rm Aut}\,\M$ acts transitively on its edges. As shown in~\cite{Jon16}, this is equivalent to the condition $\Gamma=N_{\Gamma}(M)E$, where $E$ is the Klein four-group $\langle R_0, R_2\rangle\le\Gamma$. In this case, since $N_{\Gamma}(M^+)\ge N_{\Gamma}(M)$, $\widetilde\M$ is also edge-transitive. Since $E$ acts transitively on the cosets of $N_{\Gamma}(M)$, it follows that $|\Gamma:N_{\Gamma}(M)|$ divides $4$. We therefore have

%If $\M$ is edge-transitive, in class~$T$, then the fact that $N_{\Gamma}(M^+)\ge N_{\Gamma}(M)$ means that $\M^+$ is also edge-transitive, in a class $T^+$ covered by $T$ (see Lemma~\ref{covering}), with 
\[|{\rm Aut}\,\M^+:{\rm Aut}\,\M|=2|N_{\Gamma}(M^+):N_{\Gamma}(M)|=2, 4 \;{\rm or}\; 8.\]
The most frequent value is $2$, equivalent to $\M$ being stable. However, as the following example shows, there are also unstable edge-transitive maps.

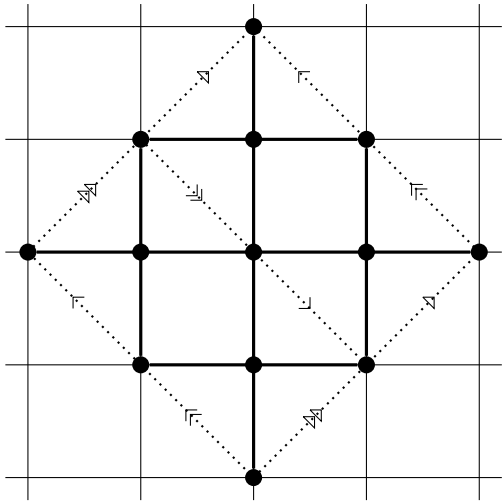
\begin{figure}[h!]
\begin{center}
\begin{tikzpicture}[scale=0.3, inner sep=0.8mm]

\draw [thin] (-11,10) to (11,10);
\draw [thin] (-11,5) to (11,5);
\draw [thin] (-11,0) to (11,0);
\draw [thin] (-11,-5) to (11,-5);
\draw [thin] (-11,-10) to (11,-10);

\draw [thin] (10,-11) to (10,11);
\draw [thin] (5,-11) to (5,11);
\draw [thin] (0,-11) to (0,11);
\draw [thin] (-5,-11) to (-5,11);
\draw [thin] (-10,-11) to (-10,11);

\node (a) at (0,10) [shape=circle, fill=black] {};
\node (b) at (-5,5) [shape=circle, fill=black] {};
\node (c) at (0,5) [shape=circle, fill=black] {};
\node (d) at (5,5) [shape=circle, fill=black] {};
\node (e) at (-10,0) [shape=circle, fill=black] {};
\node (f) at (-5,0) [shape=circle, fill=black] {};
\node (g) at (0,0) [shape=circle, fill=black] {};
\node (h) at (5,0) [shape=circle, fill=black] {};
\node (i) at (10,0) [shape=circle, fill=black] {};
\node (j) at (-5,-5) [shape=circle, fill=black] {};
\node (k) at (0,-5) [shape=circle, fill=black] {};
\node (l) at (5,-5) [shape=circle, fill=black] {};
\node (m) at (0,-10) [shape=circle, fill=black] {};

\draw [very thick] (b) to (d);
\draw [very thick] (e) to (i);
\draw [very thick] (j) to (l);

\draw [very thick] (b) to (j);
\draw [very thick] (a) to (m);
\draw [very thick] (d) to (l);

\draw [thick, dotted] (-10,0) to (0,10) to (10,0) to (0,-10) to (-10,0);
\draw [thick, dotted] (-5,5) to (5,-5);

\draw (-2.5,8) to (-2,8) to (-2,7.5) to (-2.5,8);
\draw (7.5,-2) to (8,-2) to (8,-2.5) to (7.5,-2);

\draw (-7.5,3) to (-7,3) to (-7,2.5) to (-7.5,3);
\draw (-7.8,2.7) to (-7.3,2.7) to (-7.3,2.2) to (-7.8,2.7);
\draw (2.5,-7) to (3,-7) to (3,-7.5) to (2.5,-7);
\draw (2.2,-7.3) to (2.7,-7.3) to (2.7,-7.8) to (2.2,-7.3);

\draw (2.5,8) to (2,8) to (2,7.5);
\draw (-7.5,-2) to (-8,-2) to (-8,-2.5);
\draw (2.5,-2) to (2.5,-2.5) to (2,-2.5);

\draw (7.5,3) to (7,3) to (7,2.5);
\draw (7.7,2.8) to (7.2,2.8) to (7.2,2.3);
\draw (-2.5,-7) to (-3,-7) to (-3,-7.5);
\draw (-2.3,-7.2) to (-2.8,-7.2) to (-2.8,-7.7);
\draw (-2.5,3) to (-2.5,2.5) to (-3,2.5);
\draw (-2.3,2.8) to (-2.3,2.3) to (-2.8,2.3);

\end{tikzpicture}

\end{center}
\caption{A map $\M$ with $|{\rm Aut}\,\M^+|=8|{\rm Aut}\,\M|$} 
\label{8|AutM|}
\end{figure}

\medskip

\noindent{\bf Example 2.} Let $\K$ be the regular torus map $\{4,4\}_{2,2}$ (see~\cite[\S8.3]{CM}), shown in Figure~\ref{8|AutM|} as the dotted square $|x\pm y|\le 2$ with identifications of its sides indicated by the arrows. Then ${\rm Aut}\,\K$ is the semidirect product of an abelian normal translation group of order $8$ by the stabiliser of a vertex, isomorphic to $D_4$; in particular, $|{\rm Aut}\,\K|=64$. Let $a$ be the glide reflection $(x,y)\mapsto (y+1,x+1)$ of $\K$, and let $\M=\K/\langle a\rangle$, a map of type $\{4,4\}$ on the Klein bottle represented by either of the two dotted rectangles with the indicated side-pairings. Then $\M$ is not regular, as it admits no rotations of order $4$. (There are, in fact, no regular maps on the Klein bottle, see~\cite[\S8.8]{CM}.)  

Now $\M$ admits reflections in the diagonals of each of its faces,  so it easily follows that $\M$ is edge-transitive. However, no reflection or non-identity rotation fixes the midpoint of an edge, so ${\rm Aut}\,\M$ acts regularly on the edges, with $|{\rm Aut}\,\M|=8$, the number of edges. Since $\K$ is a double cover of $\M$, and is orientable and without boundary, we have $\widetilde\M=\K$, a regular map with $|{\rm Aut}\,\widetilde\M|=8|{\rm Aut}\,\M|$. Thus $\M$ is unstable.

\iffalse
In this example, $M^+$ is the normal closure in $N(3)=\langle S_0, \ldots, S_3\mid S_i^2=1\rangle$ of the elements $(S_iS_j)^2$ for $i\ne j$, so that $N(3)/M^+\cong C_2^4$, and ${\rm Aut}\,\M^+$ is an extension of this by $\Gamma/N(3)\cong V_4$. The elements $(S_iS_j)^2$ are all in $\Gamma^+$, so $\M^+$ is orientable and without boundary. We can take $M$ to be the subgroup generated by $M^+$ and any product $S_iS_jS_k$ with distinct subscripts. Since $S_iS_jS_k$ is not in $\Gamma^+$, and does not induce a reflection of $\M^+$ (as any $S_i$ would), $\M$ is non-orientable and without boundary. We have ${\rm Aut}\,\M\cong N(3)/M\cong C_2\times C_2\times C_2$.
\fi

\medskip

\noindent{\bf Example 3.} One can extend Example~2 to an infinite series of examples by replacing $\K$ with an $m^2$-sheeted unbranched covering, namely the regular torus map $\{4,4\}_{2m,2m}$, corresponding to the subgroup $mM^+<M^+<{\mathbb Z}^2$, for any integer $m\ge 2$, and taking $a$ to be the glide reflection $(x,y)\mapsto (y+m,x+m)$. The diagram corresponding to Figure~\ref{8|AutM|} is the same, except that the square tessellation in the background is now finer, each square being subdivided into $m^2$ smaller squares. There is a similar construction using the torus map $\K=\{4,4\}_{2m,0}$ and the glide-reflection $a:(x,y)\mapsto (x+m,-y)$, but in this case the resulting unstable map $\M$ is not edge-transitive.

\medskip

\noindent{\bf Problem.} Do other surfaces provide examples of edge-transitive maps $\M$ for which $|{\rm Aut}\,\widetilde\M|=8|{\rm Aut}\,\M|$?

%%%%%%%%%%%%%

\section{A general construction}

The above examples suggest a more general construction. Let us start with an orientable map $\K$ without boundary, such that there is an orientation-reversing involution $a\in G={\rm Aut}\,\K$, and define $\M=\K/\langle a\rangle$. If $K$ is a map subgroup of $\Gamma$ for $\K$, then there is a map subgroup $M$ for $\M$ containing $K$ with index $2$, corresponding to the subgroup $\langle a\rangle\le G\cong N_{\Gamma}(M)/M$. Since $a$ reverses orientation we have $M\not\le\Gamma^+$, and since $K\le\Gamma^+$ with $|M:K|=2$ we must have $K=M^+$ so that $\K=\widetilde\M$. Conversely, every orientable double cover $\widetilde\M\to\M$ must arise in this way for some $\K$ and $a$. Now $\M$ has a non-empty boundary if and only if $a$ has a fixed point in $\K$, or equivalently, $M$ contains a conjugate of some generator $R_i$ of $\Gamma$; otherwise, $\M$ is non-orientable and without boundary, having the same type as $\K$ and having genus $g+1$ if $\K$ has genus $g$.

Here ${\rm Aut}\,\M$ is identified with $C_G(a)/\langle a\rangle$, where $C_G(a)$ is the centraliser of $a$ in $G$, so $\M$ is stable if and only if $a$ is in the centre $Z(G)$ of $G$. In order to find an unstable map without boundary, we therefore need $G$ to have a non-central orientation-reversing involution $a$, which is not a reflection.

The following example shows that there are unstable maps of each genus $g\ge 2$. (There are none on the real projective plane, since an orientation-reversing involution of a spherical map is either central or a reflection.)

\medskip

\noindent{\bf Example 4.} Let $\K$ be the regular orientable map $\{n,n\}_2$ where $n=2m$ is even (see~\cite[\S8.6 and Table~8]{CM}). This map is the Petrie dual of the spherical map $\{2,n\}$ in Example~1, obtained by replacing the faces of that map with its Petrie polygons, closed zig-zag paths turning alternately first left and first right. Then $\K$ has genus $m-1$ and automorphism group
\[G=\langle r_i\;(i=0,1,2)\mid r_i^2=(r_0r_1)^n=(r_2r_0)^2=(r_1r_2)^n=(r_0r_1r_2)^2=1\rangle\]
of order $4n$, a semidirect product of the orientation-preserving subgroup
\[G^+={\rm Aut}^+\K=\langle x=r_1r_2, y=r_2r_0\mid x^n=y^2=[x,y]=1\rangle\cong C_n\times C_2\]
by $\langle r_2\rangle\cong C_2$, with $r_2$ inverting all elements of ${\rm Aut}^+\K$ by conjugation. The involution $a=r_0r_1r_2$ reverses orientation, and it is not conjugate to any $r_i$ since these four involutions have different images in the abelianisation $C_2\times C_2\times C_2$ of $G$, so $\M=\K/\langle a\rangle$ is a non-orientable map with empty boundary and with $\widetilde\M=\K$. It has type $\{n,n\}$ and genus $m$. The centraliser of $a$ in $G^+$ is the Klein four-group $\langle x^m, y\rangle$, so $|C_G(a)|=8$ and hence $|{\rm Aut}\,\M|=4$. Since  $|{\rm Aut}\,\widetilde\M|=8m$, this shows that $\M$ is unstable for each $m\ge 2$.

\medskip

The following example shows that a single map $\mathcal K$ can have the form $\widetilde\M$ for arbitrarily many non-isomorphic maps $\M$.

\medskip

\noindent{\bf Example 5.}
Choose any integer $n\equiv 3$ mod~$(4)$ with $n\ge 11$, and define involutions $r_i\in S_n$ for $i=0, 1$ and $2$ by
\[r_0=(1,2),\quad r_1=(1)(2,n)(3,n-1)\ldots, \quad r_2=(1,2)(3,n)(4,n-1)\ldots,\]
so that $r_0$ and $r_2$ commute. Since $r_1r_2=(1, 2, \ldots, n)$ we have $\langle r_0, r_1r_2\rangle=S_n$, giving an epimorphism $\theta:\Gamma\to S_n$, $R_i\mapsto r_i$. Let $K$ be its kernel and let $\K$ be the corresponding regular map, with ${\rm Aut}\,\K\cong S_n$. Each $r_i$ is an odd permutation, so $K\le\Gamma^+$ and hence $\K$ is orientable and without boundary. Since $r_0r_1=(1,n,2)(3,n-1)(4,n-2)\ldots$, $\K$ has type $\{6,n\}$ and genus $1+\frac{(n-1)!(n-3)}{6}\sim n!/6$.

Now let $\M=\K/\langle a\rangle$ where $a=(1,2)(3,4)\ldots(m-1,m)$ for some $m\equiv 2$ mod~$(4)$ with $6\le m\le n-5$. The cycle structure of $a$ differs from that of each $r_i$, so $\M$ has empty boundary. Since $a$ is an odd permutation, $\M$ is non-orientable, with $\K=\widetilde\M$. Since $C_{S_n}(a)\cong (S_2\wr S_{m/2})\times S_{n-m}$, where $\wr$ denotes a wreath product, we have
\[|{\rm Aut}\,\M|=2^{m/2}(m/2)!(n-m)!/2\,.\]

We thus obtain $(n-7)/4$ unstable maps $\M$ for $m=6, 10, \ldots, n-5$. Any isomorphism between two of them would lift to an automorphism of $\mathcal K$, implying a conjugation between the corresponding involutions $a$. However, these all have different cycle structures, so the maps are mutually non-isomorphic.

%%%%%%%%%%%

\section{Edge-labelled graphs}

For any map $\M$, the dual of its barycentric subdivision is an embedding, in the same surface, of a permutation diagram $\mathcal D$ for the action of $\Gamma$ on $\Phi$, or equivalently of a coset diagram for $M$ in $\Gamma$: this graph has trivalent vertices corresponding bijectively to the flags of $\M$ (or cosets of $M$), and edges labelled $R_i$ indicating the action on them of the generators of $\Gamma$. It is bipartite if and only if $\M$ is orientable and without boundary; otherwise, it has a connected bipartite double, namely the corresponding permutation diagram $\widetilde{\mathcal D}$ for $\widetilde\M$.

We have ${\rm Aut_e}\,\mathcal D={\rm Aut}\,\M$ and ${\rm Aut_e}\,\widetilde{\mathcal D}={\rm Aut}\,\widetilde\M$, where the groups on the left consist of the graph automorphisms preserving the edge-labelling; this is because in either case each group consists of the permutations of $\Phi$ commuting with the action of $\Gamma$. It follows that the stability properties of maps are reflected in the corresponding properties of their associated edge-labelled graphs.

However, if we ignore the edge-labelling, the automorphism groups of the unlabelled graphs could be larger. For instance, the self-duality of the maps in Example~4 implies that there are graph automorphisms transposing the sets of edges labelled $R_0$ and $R_2$, corresponding to an outer automorphism $R_i\mapsto R_{2-i}$ of $\Gamma$ preserving the map subgroup $M$ (see~\cite{JT} for this interpretation of duality and other map operations such as Petrie duality).

Nevertheless, it is easily seen that if ${\rm Aut}\,\mathcal D={\rm Aut}\,\M$ and $\M$ is unstable, then $\mathcal D$ is also unstable, and that the converse holds if ${\rm Aut}\,\widetilde{\mathcal D}={\rm Aut}\,\widetilde\M$.

%%%%%%%%%%%%

\section{Hypermaps}

This concept of stability extends naturally to hypermaps, regarded as objects in their own right or as bipartite maps via the Walsh construction~\cite{Wal}. For an algebraic approach, we can simply omit the relation $(R_0R_2)^2=1$ from the presentation of $\Gamma$, giving a group $\Delta\cong C_2*C_2*C_2$ whose transitive permutation representations correspond to connected hypermaps~\cite{Jon16b}. The preceding ideas and constructions then carry over to hypermaps in the obvious way.

\medskip

\noindent{\bf Example 6.} If we redefine $r_0$ in Example~5 to be the involution $(2,3)$, then $r_1r_2$, $r_2r_0$ and $r_0r_1$ have orders $n$, $4$ and $4$, so $\K$ is now an orientable hypermap of type $(n,4,4)$. Factoring out the subgroups $\langle a\rangle$ of ${\rm Aut}\,\K=S_n$ as before, we obtain $(n-7)/4$ unstable hypermaps with canonical double cover $\K$.

\medskip

This concept of stability can be extended further to any permutational category~\cite{Jon16b} in which the parent group $\Gamma$ has a subgroup $\Gamma^+$ of index $2$. For instance, this is a natural extension in the theory of abstract polytopes~\cite{MS}, where one can take $\Gamma$ to be the string Coxeter group of a given Schl\"afli type, with $\Gamma^+$ its even subgroup.

%%%%%%%%%%

\section{Klein surfaces and Riemann surfaces}

The ideas and examples discussed here have close analogues in the theories of Klein surfaces and Riemann surfaces~\cite{AG, BCGG, BEGG}. If a Klein surface $\mathcal S$ is not a Riemann surface, that is, it is non-orientable or has a non-empty boundary, then it has a canonical double cover $\widetilde{\mathcal S}$ which is a Riemann surface, called its {\em complex double}. The group ${\rm Aut}\,\widetilde{\mathcal S}$ of conformal and anti-conformal automorphisms of $\widetilde{\mathcal S}$ has a subgroup ${\rm Aut}\,\mathcal S\times C_2$, which may or may not be the whole group. Special cases of these two possibilities have been considered in several papers such as~\cite{BCC, BCGS, May, PC}. It may seem impertinent to suggest a name for a concept which has been studied, apparently without a name, for several decades, but it would be consistent with the rest of this paper to call $\mathcal S$ stable or unstable as ${\rm Aut}\,\widetilde{\mathcal S}={\rm Aut}\,\mathcal S\times C_2$ or not.

%In particular, there is a faithful action of $\Gamma$ as a group of isometries of the hyperbolic plane $\mathbb H$, with the generators $R_i$ acting as reflections in a triangle with internal angles $0, \pi/2$ and $0$. The complex structure on $\mathbb H$ induces (after appropriate compactification) the structure of a Klein surface $\mathcal S=\mathbb H/M$ on any map $\M$ with map subgroup $M\le\Gamma$, and this is a Riemann surface in the case of a subgroup of $\Gamma^+$. In particular, applying this to $M^+=M\cap\Gamma^+$ shows that the surface underlying the canonical double $\widetilde\M$ of $\M$ is the complex double $\widetilde{\mathcal S}$ of $\mathcal S$.

If a map $\M$ has type $\{n,m\}$, the action of $\Gamma$ on $\Phi$ can be factored through the extended triangle group $\Delta=\Delta[m,2,n]$ obtained by adding the relations $(R_1R_2)^m=(R_0R_1)^n=1$ to the presentation of $\Gamma$. This group acts on a simply connected Riemann surface $\mathcal U$, namely the Riemann sphere, complex plane or hyperbolic plane, as a group of isometries generated by reflections in the sides of a triangle with internal angles $\pi/m$, $\pi/2$ and $\pi/n$. If $M$ is a map subgroup of $\Delta$ for $\M$, the underlying surface of $\M$ can be identified with the Klein surface $\mathcal S=\mathcal U/M$, which is a Riemann surface in the case of a subgroup of $\Delta^+=\Delta(m,2,n)$. In particular, if $\M$ is not a Riemann surface then applying this to $M^+=M\cap\Delta^+$ shows that the surface underlying the canonical double $\widetilde\M$ of $\M$ is the complex double $\widetilde{\mathcal S}$ of $\mathcal S$.

However, this analogy is not exact: we have ${\rm Aut}\,\M\le{\rm Aut}\,\mathcal S$ and  ${\rm Aut}\,\widetilde\M\le{\rm Aut}\,\widetilde{\mathcal S}$, but these inclusions could be proper. This is most obvious in cases where the underlying surfaces are the sphere, torus, real projective plane or Klein bottle (see Examples 2 and 3), since these have uncountable automorphism groups whereas the maps they support all have finite automorphism groups. Even in compact cases of higher genus, where the surfaces have finite automorphism groups, there are examples of proper inclusions, with the surface having automorphisms which do not preserve the map.

\medskip

\noindent{\bf Example 7.} It follows from the self-duality of the maps $\M$ and $\widetilde\M$ of type $\{n,n\}$ in Example~4 that their automorphism groups have index $2$ in the automorphism groups of their median maps $\M^*$ and $\widetilde\M^*$: these are maps of type $\{n,4\}$ on the same surfaces $\mathcal S$ and $\widetilde{\mathcal S}$ as $\M$ and $\widetilde\M$, but with a vertex $v_e$ at the mid-point of each edge $e$ of $\M$ or $\widetilde\M$, and edges between vertices $v_e$ and $v_{e'}$ on consecutive edges $e$ and $e'$ of a face of $\M$ or $\widetilde\M$.  

\medskip

As was the case with graphs in \S5, if ${\rm Aut}\,\widetilde{\mathcal S}={\rm Aut}\,\widetilde\M$ and $\M$ is stable, then
%${\rm Aut}\,\widetilde{\mathcal S}={\rm Aut}\,\mathcal S\times C_2$
$\mathcal S$ is stable, and if ${\rm Aut}\,\mathcal S={\rm Aut}\,\M$ and $\M$ is unstable, then
%${\rm Aut}\,\widetilde{\mathcal S}>{\rm Aut}\,\mathcal S\times C_2$
$\mathcal S$ is unstable. The latter applies in Example~5, for instance, since the maximality of $\Delta(n,2,6)$ among Fuchsian groups for $n\ge 7$  (see~\cite{Sin}) implies that ${\rm Aut}\,\mathcal S={\rm Aut}\,\M$. It also applies in Example~6, since $\M^*$ is unstable and the maximality of $\Delta(4,2,n)$ for $n\ge 5$ implies that ${\rm Aut}\,\mathcal S={\rm Aut}\,\M^*$.

\bigskip

\noindent{\bf Acknowledgement} The author thanks David Singerman for a number of helpful comments.

%%%%%%%%%%%%%%%%%%%%%
%%%%%%%%%%%%%%%%%%%%%

\end{document}